\documentclass[12pt,a4paper]{article}
\usepackage{amsfonts}
\usepackage{amsmath,amsfonts,amsthm,mathrsfs,grffile,bm}
\usepackage{epstopdf}
\usepackage{graphicx}
\usepackage{verbatim}
\usepackage[ruled]{algorithm2e}   
\usepackage{color}
\usepackage{ulem}
\usepackage{picinpar}
\usepackage{booktabs}
\usepackage[doublespacing]{setspace}
\usepackage[authoryear]{natbib}

\newtheorem{thm}{Theorem}
\newtheorem{cor}{Corollary}
\newtheorem{lem}{Lemma}

\theoremstyle{definition}   %

\theoremstyle{remark}

\numberwithin{equation}{section}

\newcommand{\ba}{\mathbf{a}}

\newcommand{\bbeta}{\boldsymbol{\beta}}

\def\eop{{\hfill\vbox{\hrule height .3pt
      \hbox{\vrule width.3pt height 7pt
      \kern 7pt
      \vrule width .3pt}
      \hrule height .3pt}} \par\bigskip}

\title{Learning rates for partially linear support vector machine in high dimensions }

\author{Yifan Xia$^{\dag}$, Yongchao Hou$^{\dag}$, Shaogao Lv$^{\ddag,*}$
\\
{\small $^{\dag}$College of Statistics,  Southwestern University of
Finance and Economics, Chengdu, China;} \\
{\small $^{\ddag}$Department of Statistics and Mathematics, Nanjing Audit University, Nanjing, China;}\\
{\small $^*$Corresponding author: kenan716@mail.ustc.edu.cn }}
\date{}

\begin{document}
\maketitle

\setcounter{page}{1}
\begin{abstract}
This paper analyzes a new regularized learning scheme for high dimensional partially linear support vector machine.
The proposed  approach consists of an empirical risk and the Lasso-type penalty for linear part, as well as the standard functional norm for nonlinear part. Here the linear kernel is used for model interpretation and feature selection, while the nonlinear kernel is adopted to enhance algorithmic flexibility. In this paper,
we develop a new technical analysis on the weighted empirical process, and establish the sharp learning rates for  the  semi-parametric estimator under the regularized conditions. Specially, our derived learning rates for semi-parametric SVM depend on not only the 
sample size and the functional complexity, but also the sparsity and the margin parameters.

\end{abstract}


\section{Introduction}
Support vector machine (SVM), originally introduced by \cite{Vapnik1995}, is well known to be a popular and powerful technique mainly due to its successful practical performances and nice theoretical foundations in machine learning. 
For supervised classification problems, SVM is based on the  margin-maximization principle endowed with a specified kernel, which is formulated by a nonlinear map from the input space to the feature space.

Denote $ \mathcal{X}$ and $\mathcal{Y}=\{-1,+1\}$ as the input space and 
corresponding output space, respectively. Let $(X,Y)\in \mathcal{X}\times \mathcal{Y}$ be 
a random vector drawn from an unknown joint distribution $\rho$ on $ \mathcal{X}\times \mathcal{Y}$. 
Suppose that all the observations $\{(Y_i,X_i)\}_{i=1}^{n}$ are available from
$\rho$.
 In empirical risk minimization, the standard $L_2$-norm SVM has the widely-used hinge loss 
plus $L_2$-norm penalty formulation. Recall that
  the empirical hinge loss function is defined by
 \begin{align*}
 \mathcal{R}_n(f)=\frac{1}{n}\sum_{i=1}^{n}  \phi_h(Y_if(X_i)),
 \end{align*}
where  the hinge loss is  $\phi_h(u)=(1-u )_+$,  with $u_+$ denoting the positive part of $u\in \mathbb{R}$. 
The standard SVM can be expressed as the following regularization problem
\begin{equation*}
\min_{f\in\mathcal{H}_K}\big\{\mathcal{R}_n(f)+
\lambda \|f\|_{K}^2\big\},
\end{equation*}
where $\lambda $ is the regularized parameter for controlling the  functional complexity of $\mathcal{H}_K$. Note that
$\mathcal{H}_K$ is referred to a reproducing kernel Hilbert space (RKHS), often specified in advance. See Section 2 for more details on RKHS. The book by \citet*{Steinwart2008} contains a good overview of SVMs and the  particularly related 
learning theory.

Among various kernel-based learning schemes including SVM, it is full of challenges how to select a suitable kernel  and there are not any perfect answers for such problem until now. See related work on kernel learning \citep*{Lancjriet2004, Micchelli2005, Wu2007, Kloft2011, Ying2009} for instances. In this paper, we consider a semi-parametric SVM problem of the linear kernel plus a general nonlinear kernel. Indeed, partial linear models in statistics have received a great attention in the last several decades, see \citep*{Muller2015, Hardle2007, Speckman1988}. Particularly, the linear part in the partial linear modelsaims at the model interpretation, and the nonlinear part is used to enhance the model flexibility. As a concrete example in stock market,   the future return of a stock $Y$ may depend on several company management indexes (e.g. shareholders structure) which are homogeneous for all the companies, and we allow linear relation with $Y$. However, the other features (e.g., from financial statements) should be nonlinear to the response, in that a company has a complex curve in terms of operation or profit pattern. In practice, load forecasting  using semi-parametric SVM gets a better prediction than the conventional way \citep*{ Abhisek2009}.  The   semi-parametric SVM are also successfully applied to analyze
pharmacokinetic and pharmacodynamic data \citep*{Seok2011}. However, to the best of our knowledge, the theoretical research on the semi-parametric support vector machines is still lacking, and this paper focuses on this topic in high dimensional setting.

High dimensional case refers to the setting where the ambient dimension $p$ of the  covariates  is very large (e.g. $p\gg n$), but only a small subset of the covariates are significantly relevant to the response. The high  dimensional estimation and inference for various models have been investigated in the last years, and the interested readers can refer to  two related book written by \cite{Buhlmann2019} and \cite{Giraud2014}. Specially, the high dimensional inference for the linear (or additive) SVM  has been wildly studied in recent years, see \citep*{Tarigan2006, Zhao2012, Zhang2016, Peng2016}.
Precisely, \cite{Tarigan2006} consider a $\ell_1$-penalized parametric estimation in high dimensions for SVM and prove the convergence rates of the excess risk term under regularity conditions. Similarly, \cite{Zhao2012} propose a group-Lasso
type regularized approachs for the nonparametric additive SVM, and provide the oracle properties of the estimator and develop an efficient numerical algorithm to compute it. For high dimensional linear SVM, \citet*{Zhang2016} and \cite{Peng2016}  explicitly investigate the statistical performance of the $\ell_1$-norm and non-convex-penalized SVM such as variable selection consistency. However, all the aforementioned works only consider a single kernel in high dimensions.
By contrast, a partially linear SVM has to consider the mutual correlation between these two kernels with different structures, and also considers the mutual effects between the sparsity and the nonlinear functional complexity.  So the non-asymptotic analysis of such semi-parametric models in high-dimensional SVM appears to be considerably more complicated than those  based on a single kernel.

Under our partial linear setting, the whole input feature  consists of two parts:
$X=(Z,T)'$, where $Z\in \mathbb{R}^p$ has a linear relation to the response, while the  sub-feature $T$ has a nonlinear effect to the response. 
Given all the observations  $\{(Y_i,Z_i,T_i)\}_{i=1}^{n}$ with the  sample size $n$, we consider a two-fold regularized learning scheme for the high dimensional PLQR,  and the semi-parametric  estimation pair $(\hat \bbeta, \hat g)$ is the unique solution by minimizing the following unconstrained optimization
\begin{equation}\label{submethod}
\min_{(f=\bbeta'Z+ g)\in\mathcal{F}}\Big\{\mathcal{R}_n(f)+\lambda_n\|\bbeta\|_1+
\mu_n \|g\|_{K}^2\Big\},
\end{equation}
where $(\lambda_n,\mu_n)
$ are two regularized  hyper-parameters for controlling the coefficients of the sparsity and  functional complexity, respectively. 
In the partial linear setting, the adopted hypothesis space $\mathcal{F}$ for SVM is a summation of the linear kernel and  the general nonlinear kernel. More precisely,
$$
\mathcal{F}:=\{f(X)=\bbeta'Z+g(T),\,\,\bbeta\in \mathbb{R}^p, g\in \mathcal{H}_K\}.
$$

To investigate the statistical performance of the proposed semi-parameter estimator \eqref{submethod}, we introduce a population target function for the partial linear SVM within $\mathcal{F}$. In this paper, 
the target function we will focus on is  a global solution $f^*$ of the following population minimization  on $\mathcal{F}$,
\begin{align}
\min_{f\in\mathcal{F}}\mathcal{R}(f),\quad \hbox{where}\,\,\mathcal{R}(f):=\mathbb{E}_\rho[\phi_h(Yf(X))].
\end{align}
Under the partial linear framework, $f^*$ can be written as: $f^*(X)=(\bbeta^*)'Z+g^*(T)$, where 
$g^*$ is the nonparametric component, belonging to a specific RKHS that will be defined in Section 2.
For the parametric part, one often assumes that the structure of $\bbeta^*$ is sparse  under high dimensional setting,
in sense that the cardinality of $S=\{j,\,\beta_j^*\neq 0,j=1,2,...,p\}$ is far less than the ambient dimension $p$.
Note that, the target function $f^*$ is quite different from the Bayes rule, and the latter  is an optimal decision function taken over all the measurable functions. We can  treat $f^*$ as a sparse approximation to the Bayes rule within $\mathcal{F}$, particularly when the true function is not sparse. In the current  literatures, we are not concerned with any approximation error induced by sparse approximation or kernel misspecification.

In this paper, we are primarily concerned with learning rates of the excess risk $\mathcal{R}(\hat f)-\mathcal{R}(f^*)$
and the estimation errors of the parametric estimator and the nonparametric estimator for the high dimensional SVM. Interestingly, the  theoretical results reveal that our derived rate of the parametric estimator depends on not only 
the sample size and the sparsity parameter, but also the functional complexity generated by the non-parametric component
and vice versa. As a byproduct, we develop a new weighted empirical process  to refine our analysis. This is one of the  key theoretical tools in the high dimensional literatures of the semi-parametric  estimation.

The rest of this article is organized as follows. In Section 2, we introduce some basic notations on RKHS that is used  to characterize the functional complexity. Then we impose some regular assumptions required to establish the convergence rates. In the end of Section 2, we explicitly propose our main theoretical results in terms of the excess risk and estimation errors. Section 3 is devoted to a detailed proof for the main theorems, and also proves some useful lemmas associated with the weighted empirical process. Section 4 concludes this paper with discussions and future possible researches.

{\bf Notations}. 
We use $[p]$ to denote the set $\{1,2,...,p\}$. For a vector $\ba=(a_1,a_2,...,a_p)\in \mathbb{R}^p$,  the $\ell_q$-norm is defined as $\|\ba\|_q=\big(\sum_{i\in[p]}|a_i|^q\big)^{1/q}$.   For two sequences of numbers $a_n$  and  $b_n$, we use an $a_n=O(b_n)$ to denote that $a_n\leq Cb_n$ for some finite positive constant $C$ for all $n$. If 
both $a_n=O(b_n)$ and $b_n=O(a_n)$, we use the notation $a_n\simeq b_n$. We also use an $a_n=\Omega(b_n)$ for $a_n\geq Cb_n$.
For a function, we 
denote the $L_2$-norm of $f$ by $\|f\|_2=\big(\int_X f(x)^2d\rho_X(x)\big)^{1/2}$ with some distribution $\rho_X$.

\section{Conditions and Main Theorems}
We begin with the background and notation required for the main statements of our problem.
First of all, we introduce the notation of RKHS. RKHS can be defined by any symmetric and positive semidefinite kernel function $K:\,\mathcal{T} \times \mathcal{T}\rightarrow \mathbb{R}$. For each $t \in \mathcal{T}$, the function $t'\rightarrow K(t',t)$ is contained with the Hilbert space $\mathcal{H}_K$; moreover, the Hilbert space is endowed with an inner product $\langle\cdot,\cdot\rangle_K$ such that $K(\cdot,t)$ acts as the representer of the evaluation. Especially, 
the reproducing property of RKHS plays an important role in the theoretical analysis and numerical optimization for any kernel-based method,
\begin{align}
\label{C1}
f(t)=\langle f,K(\cdot,t)\rangle_K,\quad \forall\, t \in \mathcal{T}.
\end{align}
This property also implies that $\|f\|_\infty\leq \kappa \|f\|_K$ with $\kappa:=\max_{t\in \mathcal{T}}|K(t,t)|<\infty$.  Moreover, by Mercer's theorem, a kernel $K$ defined on a compact subset of $\mathcal{T}$ admits the following eigen-decomposition,
\begin{align}
\label{C2}
K(t,t')=\sum_{\ell=1}^\infty \mu_\ell \phi_\ell(t) \phi_\ell(t'),\,\,t,t'\in \mathcal{T},
\end{align}
where $\mu_1\geq \mu_2\geq \cdots>0$ are the eigenvalues and $\{\phi_\ell\}_{\ell=1}^\infty$ is an orthonormal basis in $L_2(\rho_T)$. The decay rate of $\mu_\ell$ completely characterizes the complexity of RKHS induced by a kernel $K$, and generally it has equivalent relationships with various entropy numbers, see \citet*{Steinwart2008}
for details. With these preparations, we define the quantity,
\begin{align}
\label{C3}
\mathcal{Q}_n(r)=\frac{1}{\sqrt{n}}\Big[\sum_{\ell=1}^{\infty}\min\{r^2,\mu_\ell\}\Big]^{1/2},\quad \forall \,r>0.
\end{align}
Let $\nu_n$ be the smallest positive solution to the inequality,
$
40\nu_n^2\geq\mathcal{Q}_n(\nu_n),
$
where $40$ is only a technical constant.

Then, due to the mutual effects between the high dimensional parametric component and the nonparametric one, we  introduce the following quantity related to the convergence rates of the semi-parametric estimate, as illustrated in (\ref{C4}),
\begin{align}
\label{C4}
\gamma_n:=\max\Big\{\nu_n,\sqrt{\frac{\log p}{n}} \Big\}.
\end{align}

We now describe our main assumptions. Our first assumption deals with the tail behavior of the covariate of the linear part.

{\bf Assumption A}. (i) For simplicity,  we assume that $\|Z\|_\infty\leq C_0<\infty$ with some positive constant $C_0$;
(ii) The largest eigenvalue of $\mathbb{E}[ZZ']$ is finite, denoted by $\Lambda_{\max}>0$. 

It appears that a bound on the $Z$-values is a restrictive assumption, ruling out the standard sub-gaussian covariates.
However, we can usually approximate a non-bounded distribution with its truncated version. Imposing such assumption is only for technical simplicity and may be relaxed to general thin-tail random variables. Assumption A(ii) is fairly standard in the literature to identify the coefficients associated to $Z$.

{\bf Assumption B}. There exist the constants $C_1>0$ and $\zeta\geq 2$ such that, for all $f\in \mathcal{F}$, the equation (\ref{C5}) holds,
\begin{align}
\label{C5}
\mathcal{R}( \bbeta,  g)-\mathcal{R}(\bbeta^*, g^*)\geq C_1\|f-f^*\|_2^{\zeta}.
\end{align}
The parameter $\zeta$ is called the Bernstein parameter introduced by \citet*{Bartlett2006,Pierre2019}.  Fast rates will usually be derived when $\zeta=2$.  This condition is essentially a qualification of the identifiability condition of the objective function at its minimum $f^*$.
Note that, the Bernstein parameter is slightly different from the classical margin parameter adopted by \cite{Tarigan2006,Chen2004}.

To estimate the parametric and nonparametric parts respectively, 
 some conditions concerning correlations between $Z$ and $T$ are required.
For each $j\in[p]$, let $\Pi^{(j)}_T$ be the projection of $Z^{(j)}$ onto $\mathcal{H}_K$. To be precise, $\Pi^{(j)}_T=g_j^*(T)$ with (\ref{C6}),
\begin{align}
\label{C6}
g_j^*=\arg\min_{g\in\mathcal{H}_K}\mathbb{E}_{Z^{(j)},T}[(Z^{(j)}-g(T))^2].
\end{align}
Let $\Pi_{Z|T}:=(\Pi^{(1)}_T,...,\Pi^{(p)}_T)'$ and $Z_T=Z-\Pi_{Z|T}$. Each function $g^*_j$ can be viewed as the best approximation  of $\mathbb{E}[Z^{(j)}|T]$ within  $\mathcal{H}_K$.
In the extreme case ($Z$ is uncorrelated with $T$),
$\Pi_{Z|T}=0$. The following condition is quite common in the semi-parametric estimation \citep*{Muller2015}, ensuring that there is enough information in the data to identify the parametric coefficients.

{\bf Assumption C}.   The smallest eigenvalue  of
$\mathbb{E}[Z_TZ_T']$ is  bounded below by a constant $\Lambda_{\min}>0$.

Note that, the equation (\ref{projection}) always holds with the definition of  projection on the  $\|\cdot\|_2$-norm,
\begin{align}\label{projection}
\|\bbeta'Z+g(T)\|_2^2=\|\bbeta'Z_T\|_2^2+|\bbeta'\Pi_{Z|T}+g(T)\|_2^2,\quad \forall \,f\in \mathcal{F}.
\end{align}
This equality ensures that the parametric estimation can be separated from the total estimation, which is very useful in our proof.

We are in a position to derive the learning rate of the estimator $(\hat \bbeta,\hat g)$ defined by minimization \eqref{submethod}. We  allow that the number of dimension $p$ and the number of active covariates $s:=|S|$ which are increasing with respect to the sample size $n$, while $s\ll p$ and the dimension of $T$  is fixed.

\begin{thm}
Let  $(\hat \bbeta,\hat g)$ be the proposed semi-parametric estimator for SVM defined in \eqref{submethod}, with the regularization parameters  $\lambda_n=\sqrt{\log p/n}$ and $\mu_n\simeq \gamma_n^2$. If Assumptions A, B, and C hold, the equation \eqref{C7} holds with the probability at least $1-2p^{-A/2}-16p^{3-A}$ with some $A>3$, 
\begin{align}
\label{C7}
\mathcal{R}(\hat \bbeta, \hat g)-\mathcal{R}(\bbeta^*, g^*)=O\big((\gamma_n+\sqrt{s\log p/ n})^{\frac{\zeta}{\zeta-1}}\big),
\end{align}
and at the meantime the estimation error has the form \eqref{C8},
\begin{align}
\label{C8}
\|\hat \bbeta-\bbeta^*\|_2=O\big((\gamma_n+\sqrt{s\log p/ n})^{\frac{1}{\zeta-1}}\big),\quad 
\|\hat g-g^*\|_2
=O\big((\gamma_n+\sqrt{s\log p/ n})^{\frac{1}{\zeta-1}}\big).
\end{align}
\end{thm}

Remark that, this rate may be interpreted as the sum of a subset selection term ($\sqrt{s\log p/ n}$) for the linear part and a fixed dimensional non-parametric estimation term ($\nu_n$). Depending on the scaling of the triple $(n,p,s)$ and the smoothness of the  RKHS $\mathcal{H}_K$, either the subset selection term or the non-parametric estimation term may dominate the estimation. In general, if $s \log p/ n = o(\nu_n^2)$, the $s$-dimensional parametric term can dominate the estimation, so can the vice versa otherwise. At the boundary, the scalings of the two terms are equivalent. In the best situation ($\zeta=2$), our derived rate of the excess risk is the same as the optimal rate achieved by those least square approaches, see \citep*{Koltchinskii2010,Muller2015}
for details.

Note also that, it is easy to check that Theorems 1 still holds if $p$ in the confidence probability is replaced by an arbitrary $\tilde{p}\geq p$ such that $\log \tilde{p} \geq 2\log\log n$. In this case, the divergence of $p$ is not needed and the probability bounds in the theorem becomes $1-2\tilde{p}^{-A/2}-16\tilde{p}^{3-A}$.

A number of corollaries of Theorem 1 can be obtained with particular choices of different kernels. First of all, we present finite-dimensional  $m$-rank operators, i.e., the kernel function $K$ can be expressed in terms of $m$ eigenfunctions. These eigenfunctions include the linear functions, polynomial functions, as well as the function class based on finite dictionary expansions.

\begin{cor}
Under the same conditions as Theorem 1, consider a nonlinear kernel with finite rank $m$. Then the semi-parametric estimator  for SVM defined in \eqref{submethod} with $\lambda_n=\sqrt{\log p/n}$ and $\mu_n\simeq \gamma_n^2$ satisfies the condition \eqref{C9},
\begin{align}
\label{C9}
\mathcal{R}(\hat \bbeta, \hat g)-\mathcal{R}(\bbeta^*, g^*)=O_p\Big(\Big(\frac{s\log p}{n}+\frac{m}{n}\Big)^{\frac{\zeta}{2(\zeta-1)}}\Big),
\end{align}
where
$$
\|\hat \bbeta-\bbeta^*\|_2=O_p\Big(\Big(\frac{s\log p}{n}+\frac{m}{n}\Big)^{\frac{1}{2(\zeta-1)}}\Big),\quad 
\|\hat g-g^*\|_2
=O_p\Big(\Big(\frac{s\log p}{n}+\frac{m}{n}\Big)^{\frac{1}{2(\zeta-1)}}\Big).
$$
\end{cor}
For a finite rank kernel and for any $r>0$, we have
$\mathcal{Q}_n(r)\leq r\sqrt{\frac{m}{n}}$, which follows by the result of Theorem 1. Corollary 1
corresponds to the linear case for SVM when $s\simeq m$.
The existing theory in the literatures on the linear SVM has paid constant attention to the analysis of the generalization error and 
variable selection consistency.
\cite{Zhang2016} considers the non-convex penalized SVM in terms of the variable selection consistency and oracle property in high dimension, however, their results are based on a restrictive condition in case of $p\ll n^{1/2}$. So the ultra-high dimensional cases ($p=O(e^{n^r})$ with $r<1$) are excluded. Under the constrained eigenvalues constant condition, \cite{Peng2016} provides a tight upper bound of the linear SVM estimator in the $\ell_2$ norm, with an order $\sqrt{s\log p/n}$, which is the same as our rate in Corollary 1 when $\zeta=2$. 

Secondly, we state a result for the RKHS with  countable eigenvalues,  decaying at a rate $\mu_\ell \simeq(1/\ell)^{2\alpha}$ for some smooth parameter $\alpha > 1/2$. In fact, this type of scaling covers the Sobolev spaces,  consisting of derivative functions with $\alpha$.

\begin{cor}
Under the same conditions as Theorem 1, consider a kernel with the eigenvalue decay $\mu_\ell \simeq(1/\ell)^{2\alpha}$ for some  $\alpha > 1/2$. Then the semi-parametric estimator defined in \eqref{submethod} with $\lambda_n=\sqrt{\log p/n}$ and $\mu_n\simeq \gamma_n^2$ satisfies the equation \eqref{C10},
\begin{align}
\label{C10}
\mathcal{R}(\hat \bbeta, \hat g)-\mathcal{R}(\bbeta^*, g^*)=O_p\Big(\Big(\frac{s\log p}{n}+n^{-\frac{2\alpha}{2\alpha+1}}\Big)^{\frac{\zeta}{2(\zeta-1)}}\Big),
\end{align}
where
$$
\|\hat \bbeta-\bbeta^*\|_2=O_p\Big(\Big(\frac{s\log p}{n}+n^{-\frac{2\alpha}{2\alpha+1}}\Big)^{\frac{1}{2(\zeta-1)}}\Big),\quad 
\|\hat g-g^*\|_2
=O_p\Big(\Big(\frac{s\log p}{n}+n^{-\frac{2\alpha}{2\alpha+1}}\Big)^{\frac{1}{2(\zeta-1)}}\Big).
$$
\end{cor}
In the previous corollary, we need to compute the critical univariate rate $\nu_n$. Given the assumption of polynomial eigenvalue decay, a truncation argument shows that $\mathcal{Q}_n(r)=O(\frac{r^{\frac{2\alpha}{2\alpha-1}}}{\sqrt{n}})$, i.e., $\nu_n^2\simeq n^{-\frac{2\alpha}{2\alpha+1}}$. As opposed to Corollary 1, 
we now discuss the special case where the functional complexity dominates the esimation, that is, the rate of the excess risk is an order $O\big(n^{-\frac{\alpha}{2\alpha+1}}\big)^{\frac{\zeta}{\zeta-1}}$. This is a better rate campared with those in \cite{Chen2004} and  \cite{Wu2007}. The learning rate in \cite{Chen2004} is
derived as $O_p(n^{\frac{\theta}{\theta+\eta+\theta \eta}})$, where $\theta$ is a separation parameter corresponding to $\zeta=\frac{1+2\theta}{\theta}$, and $\eta$ is a power appearing in the covering number, satisfying $\eta\simeq \frac{1}{\alpha}$.  Our rate can be proved to be better than that of Chen in the best case with $\zeta=2$. The similar arguments also hold when considering the result in \cite{Wu2007}.

\section{Proofs}
In this section, we provide the proofs of our  main theorem (Theorem 1). At a high-level, Theorem 1 is based on an appropriate adaptation to the semi-parametric settings of various techniques, developed for sparse linear regression or additive non-parametric estimation in high dimensions \citep*{Buhlmann2019}. In contrast to the parametric setting or additive setting, it involves structural deals from the empirical process theory to control the error terms in the semi-parametric case . In particular, we make use of several concentration theorems for the empirical processes \citep*{Geer2000}, as well as the results on the Rademacher complexity of kernel classes \citep*{Bartlett2005}. 

\subsection{Proof for Theorem 1}
We write the total empirical objective as the equation \eqref{P1},
\begin{align}
\label{P1}
\mathcal{L}(\bbeta,g)=\mathcal{R}_n(\bbeta,g)+\lambda_n\|\bbeta\|_1
+\mu_n
\|g\|_{K}^2.
\end{align}

The population risk for partial linear SVM is defined by \eqref{P2},
\begin{align}
\label{P2}
\mathcal{R}(\bbeta,g)=\mathbb{E}\big[ \phi_h(Y[\bbeta'Z+g(T)])\big].
\end{align}

According to the definition of $(\hat \bbeta, \hat g)$, it holds that
	$\mathcal{L}(\hat \bbeta, \hat g)\leq \mathcal{L}(\bbeta^*, g^*)$. That means,
	\begin{align}
\label{P3}
\mathcal{R}_n(\hat \bbeta, \hat g)+\lambda_n\|\hat \bbeta\|_1+\mu_n
	\|\hat g\|_{K}^2\leq \mathcal{R}_n(\bbeta^*, g^*)+\lambda_n\|\bbeta^*\|_1+\mu_n
	\| g^*\|_{K}^2.
	\end{align}
The inequality \eqref{P3} can be rewritten into the form \eqref{firststep},
\begin{align}\label{firststep}
&\mathcal{R}(\hat \bbeta, \hat g)-\mathcal{R}(\bbeta^*, g^*)+\lambda_n\|\hat \bbeta-\bbeta^*\|_1+\mu_n/2	\|\hat g-g^*\|_{K}^2\nonumber\\
&\leq \mathcal{R}(\hat \bbeta, \hat g)-\mathcal{R}_n(\hat \bbeta, \hat g)+\mathcal{R}_n(\bbeta^*, g^*)-\mathcal{R}(\bbeta^*, g^*)
+2\lambda_n\|(\hat \bbeta-\bbeta^*)_S\|_1+2\mu_n
\| g^*\|_{K}^2.
\end{align}
For simplicity, we denote,
\begin{align}
\label{P4}
\nu_n(f):=\nu_n(\bbeta,g)=(\mathcal{R}_n( \bbeta,  g)-\mathcal{R}( \bbeta, g)),\quad \forall \,f(X)=\bbeta'Z+g(T).
\end{align}

In order to derive the upper bound of $|\nu_n(\hat f)-\nu_n(f^*)|$ in  \eqref{P4}, a new weighted empirical process is proposed in our semi-parametric high dimensional setting. The process is relevant to the uniform law of large number in a mixed function space. The \textbf{Lemma \ref{wieghtineq}} can be derived via the peeling device which is often used in probabilistic theory.

\begin{lem}\label{wieghtineq}
	Let $\mathcal{E}$ be the event
	\begin{align}
\label{P5}
\mathcal{E}:=\Big\{|\nu_n(f)-\nu_n(f^*)|\leq D_0\Big(\sqrt{\frac{\log p}{n}}\|\bbeta-\bbeta^*\|_1+\gamma_n\| g-g^*\|_{2}+	\gamma_n^2\| g-g^*\|_{K}+e^{-p}\Big)\Big\},
\end{align}
where $D_0$ is a constant in the proof of the \textbf{Lemma \ref{wieghtineq}}.
Because $p=\Omega(\log n)$ and $p=o(e^n)$, the inequality \eqref{P6} holds for  some universal constant $A>3$,
\begin{align}
\label{P6}
\mathbb{P}(\mathcal{E})\geq 1-2p^{-A/2}-16p^{3-A}.
\end{align}
\end{lem}

We continue our proof  \eqref{firststep} along with the results established in Lemma  \ref{wieghtineq}.
Apply the weighted empirical process and we can obtain \eqref{secdstep} from \eqref{firststep},
\begin{align}\label{secdstep}
&\mathcal{R}(\hat \bbeta, \hat g)-\mathcal{R}(\bbeta^*, g^*)+\lambda_n\|\hat \bbeta-\bbeta^*\|_1+\frac{\mu_n}{2}	\|\hat g-g^*\|_{K}^2\nonumber\\
&\leq  D_0\big(\sqrt{\log p/n}\|\bbeta-\bbeta^*\|_1+\gamma_n\| g-g^*\|_{2}+	\gamma_n^2\| g-g^*\|_{K}+e^{-p}\big)\nonumber\\
&+2\lambda_n\|(\hat \bbeta-\bbeta^*)_S\|_1+2\mu_n
\| g^*\|_{K}^2.
\end{align}
Therefore, when the conditions $\lambda_n\geq 2D_0\sqrt{\log p/n}$ and $\mu_n\geq 2D_0\gamma_n^2$ are both satisfied, \eqref{thirdstep} holds after the inequality \eqref{secdstep},
\begin{align}\label{thirdstep}
&\mathcal{R}(\hat \bbeta, \hat g)-\mathcal{R}(\bbeta^*, g^*)+\frac{\lambda_n}{2}\|\hat \bbeta-\bbeta^*\|_1+\frac{\mu_n}{4}	\|\hat g-g^*\|_{K}^2\nonumber\\
&\leq  D_0\big(\gamma_n\| g-g^*\|_{2}+	\gamma_n^2/2\big)
+2\lambda_n\|(\hat \bbeta-\bbeta^*)_S\|_1+2\mu_n
\| g^*\|_{K}^2,
\end{align}
where we use the basic inequality $2uv\leq u^2+v^2$. Since $p=\Omega(\log n)$ implies that
$e^{-p}=O(n^{-1})=O(\gamma_n^2)$, \eqref{storng} can be derived with Assumption B, C and the  equality \eqref{projection},
\begin{align}\label{storng}
\mathcal{R}(\hat \bbeta, \hat g)-\mathcal{R}(\bbeta^*, g^*)\geq C_1\|\hat f-f^*\|_2^{\zeta}
\geq C_1\|(\hat \bbeta-\bbeta^*)'Z_T\|_2^{\zeta}\geq C_1\Lambda_{\min}^{\zeta/2}\|\hat \bbeta-\bbeta^*\|_2^{\zeta}.
\end{align}
Moreover, \eqref{P7} follows by Assumption A after some simple computations,  
\begin{align}
\label{P7}
\left\{
\begin{array}{l}
\| g-g^*\|_{2}\leq \| f-f^*\|_{2}+\Lambda_{\max}\|\hat \bbeta-\bbeta^*\|_2,\\
\lambda_n\|(\hat \bbeta-\bbeta^*)_S\|_1\leq \lambda_n\sqrt{s}\|(\hat \bbeta-\bbeta^*)_S\|_2\leq \lambda_n\sqrt{s}\|\hat \bbeta-\bbeta^*\|_2,
\end{array}
\right.
\end{align}
where the last inequality follows from the Cauchy-Schwartz inequality.
Substitute \eqref{P7} into \eqref{secdstep} and we can obtain the inequality \eqref{imdestep},
\begin{align}\label{imdestep}
&C_1\|\hat f-f^*\|_2^{\zeta}+\lambda_n\|\hat \bbeta-\bbeta^*\|_1+\mu_n/2	\|\hat g-g^*\|_{K}^2\nonumber\\
&\leq  D_0\gamma_n\| f-f^*\|_{2}+ D_0\gamma_n\Lambda_{\max}\|\hat \bbeta-\bbeta^*\|_2+	D_0\gamma_n^2/2
+2\sqrt{s}\lambda_n\|\hat \bbeta-\bbeta^*\|_2+2\mu_n
\| g^*\|_{K}^2.
\end{align}
 For any $\theta>0$, we can then derive the inequality \eqref{yongy} with the Young inequality ($uv\leq u^{q}/q+v^p/p$ with $1/q+1/p=1$), 
\begin{align}\label{yongy}
&\gamma_n\| f-f^*\|_{2}\leq \big(\frac{\zeta-1}{\zeta}\big)\theta^{\frac{\zeta}{1-\zeta}}\gamma_n^{\frac{\zeta}{\zeta-1}}+\frac{\theta^\zeta}{\zeta}\|\hat f-f^*\|_2^{\zeta},\\
&  c_n\|\hat \bbeta-\bbeta^*\|_2\leq \Big(\frac{\zeta-1}{\zeta}\Big)\theta^{\frac{\zeta}{1-\zeta}}c_n^{\frac{\zeta}{\zeta-1}}+\frac{\theta^\zeta}{\zeta}\|\hat \bbeta-\bbeta^*\|_2^{\zeta},
\end{align}
where $c_n:=D_0\gamma_n\Lambda_{\max}+2\sqrt{s}\lambda_n$.  \eqref{iiimdestep} holds if $\theta$ is small enough to satisfy $\theta^{\zeta}\leq \frac{C_1\zeta}{2}$,
\begin{align}\label{iiimdestep}
&\frac{C_1}{2}\|\hat f-f^*\|_2^{\zeta}+\lambda_n\|\hat \bbeta-\bbeta^*\|_1+\frac{\mu_n}{2}	\|\hat g-g^*\|_{K}^2\nonumber\\
&\leq D_0 \big(\frac{\zeta-1}{\zeta}\big)\theta^{\frac{\zeta}{1-\zeta}}\gamma_n^{\frac{\zeta}{\zeta-1}}+ 	\frac{D_0\gamma_n^2}{2}
+c_n\|\hat \bbeta-\bbeta^*\|_2+2\mu_n
\| g^*\|_{K}^2.
\end{align}
Furthermore, combine  \eqref{storng}, \eqref{yongy} with \eqref{iiimdestep} and we can conclude the inequality \eqref{P8},
\begin{align}
\label{P8}
&\frac{C_1\Lambda_{\min}^{\frac{\zeta}{2}}}{2}\|\hat \bbeta-\bbeta^*\|_2^{\zeta}+\frac{\lambda_n}{2}\|\hat \bbeta-\bbeta^*\|_1+\frac{\mu_n}{4}	\|\hat g-g^*\|_{K}^2\nonumber\\
&\leq D_0 \big(\frac{\zeta-1}{\zeta}\big)\theta^{\frac{\zeta}{1-\zeta}}\gamma_n^{\frac{\zeta}{\zeta-1}}+ 	\frac{D_0\gamma_n^2}{2}+2\mu_n
\| g^*\|_{K}^2+\big(\frac{\zeta-1}{\zeta}\big)\theta^{\frac{\zeta}{1-\zeta}}c_n^{\frac{\zeta}{\zeta-1}},
\end{align}
where the condition $\theta^{\zeta}\leq \frac{C_1\zeta\Lambda^{\zeta/2}_{\min}}{4}$ is additionnally required so that
$\frac{\theta^\zeta}{\zeta}$ is ignorable. In this case, we can derive \eqref{P9},
\begin{align}
\label{P9}
&\|\hat \bbeta-\bbeta^*\|_2^{\zeta}=O_p\big((\gamma_n+\sqrt{s}\lambda_n)^{\frac{\zeta}{\zeta-1}}+\mu_n\big),\\
&\lambda_n\|\hat \bbeta-\bbeta^*\|_1=O_p\big((\gamma_n+\sqrt{s}\lambda_n)^{\frac{\zeta}{\zeta-1}}+\mu_n\big).
\end{align}
Moreover, we will obtain \eqref{P10} based on \eqref{P9} and \eqref{iiimdestep},
\begin{align}
\label{P10}
\|\hat f-f^*\|_2^{\zeta}=O_p\big((\gamma_n+\sqrt{s}\lambda_n)^{\frac{\zeta}{\zeta-1}}+\mu_n+\gamma_n\mu_n^{\frac{1}{\zeta}}\big)
=O_p\big(\gamma_n^{\frac{\zeta}{\zeta-1}}+\gamma_n^{\frac{\zeta+2}{\zeta}}\big),
\end{align}
where we choose $\mu_n\simeq\gamma_n^2$ and $\gamma_n=\Omega(\lambda_n)$. Therefore, it is concluded that \eqref{P11} holds by the triangle inequality and Assumption C,
\begin{align}
\label{P11}
\|\hat g-g^*\|_2
=O_p\big((\gamma_n+\sqrt{s}\lambda_n)^{\frac{1}{\zeta-1}}\big).
\end{align}
Finally, plugging the  derived upper bounds into \eqref{thirdstep},
we obtain the desired upper bound of the excess risk $\mathcal{R}(\hat \bbeta, \hat g)-\mathcal{R}(\bbeta^*, g^*)$.
This completes the proof. $\Box$

\subsection{The semi-parametric weighted emprical process}
In order to prove Lemma 1, some auxiliary results is required which is on uniform law of large number or concerntation inequalities, stated as Lemmas \ref{consentration} \citep*{Massart2000}.

\begin{lem}\label{consentration}
	Let $U_1,...,U_n$ be independent and identically distributed copies of a random variable $U\in \mathcal{U}$.
	Let $\Gamma$ be a class of real-valued functions on $\mathcal{U}$ satisfying $\sup_{u}|\gamma(u)|\leq D$ for all
	$\gamma\in \Gamma$. Define
	$$
	\mathbf{Z}:=\sup_{\gamma\in \Gamma}\Big|\frac{1}{n}\sum_{i=1}^n\{\gamma(U_i)-\mathbb{E}[\gamma(U_i)]\}\Big|,
	$$ 
	and
	$$
B^2:=\sup_{\gamma\in \Gamma}var(\gamma(U)).
	$$
	Then there exists a universal constant $N_0$ such that 
	$$
	\mathbb{P}\Big(\mathbf{Z}\geq N_0 \big[ \mathbb{E}[\mathbf{Z}]+B\sqrt{r/n}+D r/n\big]\Big)\leq \exp(-r),\quad \forall\,r>0. 
	$$
\end{lem}

{\bf Proof for Lemma 1.}
For any $f(x)= \bbeta'z+g(t)$, we define \eqref{P12} to apply Lemma \ref{consentration},
\begin{align}
\label{P12}
\gamma(u)=\phi_h(y[\bbeta'z+g(t)])-\phi_h(y[(\bbeta^*)'z+g^*(t)]),\,\,\bbeta\in \mathbb{R}^p \,\hbox{and} \,g\in \mathcal{H}_K.
\end{align}
Based on \eqref{P12}, a bounded set of functions is introduced,
$$
\Gamma_\Delta:=\big\{\gamma,\,\,\sqrt{\log p/n}\|\bbeta-\bbeta^*\|_1\leq \Delta_\beta, \gamma_n\| g-g^*\|_{2}\leq \Delta_{-}, \gamma_n^2\| g-g^*\|_{K}\leq \Delta_{+}\big\},
$$
where we write the triplet $\Delta=( \Delta_\beta,\Delta_{-},\Delta_{+})$ .
Since $\phi_h(\cdot)$ in \eqref{P12} is Lipschitz with constant $1$, the inequality \eqref{P13} holds for any $u$,
\begin{align}
\label{P13}
\begin{array}{rl}
|\gamma(u)|&\leq |(\bbeta-\bbeta^*)'z|+|g(t)-g^*(t)|\leq C_0\|\bbeta-\bbeta^*\|_1+\kappa \| g-g^*\|_{K}\\
&\leq C_0\sqrt{\frac{n}{\log p}}\Delta_\beta+\kappa\frac{n}{\log p}\Delta_{+}.
\end{array}
\end{align}
\eqref{P13} implies that  if we take \eqref{dbound} in Lemma \ref{consentration},
\begin{align}\label{dbound}
D:=C_0\sqrt{\frac{n}{\log p}}\Delta_\beta+\kappa\frac{n}{\log p}\Delta_{+},
\end{align}
and \eqref{P14} is also derived by the Lipschitz property,
\begin{align}
\label{P14}
\begin{array}{rl}
B^2&\leq 2\mathbb{E}[((\bbeta-\bbeta^*)'Z)^2]+2\mathbb{E}[(g(T)-g^*(T))^2]\\
&\leq \frac{2n}{\log p}\big(C_0^2\Delta_{\beta}^2+\Delta_{-}^2\big),
\end{array}
\end{align}
we can plug \eqref{dbound} and \eqref{boubund} into Lemma \ref{consentration} to yield \eqref{immderuslt},
\begin{align}\label{boubund}
B:=\sqrt{\frac{2n}{\log p}}\big(C_0\Delta_{\beta}+\Delta_{-}\big),
\end{align}
\begin{align}\label{immderuslt}
\mathbb{P}\Big(\mathbf{Z}\geq N_0 \Big[ \mathbb{E}[\mathbf{Z}]+\sqrt{\frac{2r}{\log p}}\big(C_0\Delta_{\beta}+\Delta_{-}\big)+\frac{\kappa r}{\log p}\Delta_{+}\Big]\Big)\leq \exp(-r),\quad \forall\,r\in(0,2n). 
\end{align}

It remains to provide the upper bound of the term $\mathbb{E}[\mathbf{Z}]$.
Let $\sigma_1,...,\sigma_n$ be a Rademacher sequence independent of $(Y_1,X_1),...,(Y_n,X_n)$.  The inequality \eqref{P15} can be obtained by symmetrization and
the contraction inequality,
\begin{align}
\label{P15}
\begin{array}{rl}
\mathbb{E}[\mathbf{Z}]&\leq 4\mathbb{E}\Big(\sup_{\gamma\in \Gamma_\Delta}\Big|\frac{1}{n}\sum_{i=1}^n\sigma_i(f(X_i)-f^*(X_i))\Big|\Big)\\
&\leq 4\mathbb{E}\Big(\sup_{\gamma\in \Gamma_\Delta }\Big|\frac{1}{n}\sum_{i=1}^n\sigma_i((\bbeta-\bbeta^*)'Z_i)\Big|\Big)\\
&+4\mathbb{E}\Big(\sup_{\gamma\in \Gamma_\Delta}\Big|\frac{1}{n}\sum_{i=1}^n\sigma_i(g(T_i)-g^*(T_i))\Big|\Big).
\end{array}
\end{align} 
By Bernstein inequality and the union bound, we can get the inequality \eqref{uuperexp},
\begin{align}\label{uuperexp}
\mathbb{E}\Big(\sup_{\gamma\in \Gamma_\Delta }\Big|\frac{1}{n}\sum_{i=1}^n\sigma_i((\bbeta-\bbeta^*)'Z_i)\Big|\Big)
&\leq \mathbb{E}\Big(\sup_{j\in[p]}\Big|\frac{1}{n}\sum_{i=1}^n\sigma_iZ_{ij}\Big|\Big)
\sup_{\gamma\in \Gamma_\Delta }\|\bbeta-\bbeta^*\|_1\nonumber\\
&\leq C_0\Delta_{\beta}.
\end{align}
Moreover, applying Talagrand's concentration inequality once again, we get \eqref{P16} with the probability at least $1-e^{-r}$
\begin{align}
\label{P16}
\begin{array}{rl}
\mathbb{E}\Big(\sup_{\gamma\in \Gamma_\Delta}\Big|\frac{1}{n}\sum_{i=1}^n\sigma_i(g(T_i)-g^*(T_i))\Big|\Big)
&\leq
N_0\Big(\sup_{\gamma\in \Gamma_\Delta}\Big|\frac{1}{n}\sum_{i=1}^n\sigma_i(g(T_i)-g^*(T_i))\Big|\\
&+\sqrt{\frac{2r}{\log p}}\Delta_{-}+\frac{\kappa r}{\log p}\Delta_{+}\Big).
\end{array}
\end{align}
Besides the result in \citet*{Koltchinskii2010} that,
$$
\sup_{\gamma\in \Gamma_\Delta}\Big|\frac{1}{n}\sum_{i=1}^n\sigma_i(g(T_i)-g^*(T_i))\Big|\leq \nu_n\| g-g^*\|_{2}+
\nu_n^2\| g-g^*\|_{K},\quad \forall\,g\in \mathcal{H}_K,
$$
the inequality \eqref{nonlexp} holds with the probability at least $1-2e^{-r}-p^{-A/2}$, 
\begin{align}\label{nonlexp}
\mathbb{E}\Big(\sup_{\gamma\in \Gamma_\Delta}\Big|\frac{1}{n}\sum_{i=1}^n\sigma_i(g(T_i)-g^*(T_i))\Big|\Big)
&\leq
N_0\Big(\Delta_{-}+\Delta_{+}+\sqrt{\frac{2r}{\log p}}\Delta_{-}+\frac{\kappa r}{\log p}\Delta_{+}\Big),
\end{align}
where $A>3$ is some constant in Section 5 in \citet*{Koltchinskii2010}.
Thus, combining \eqref{uuperexp}, \eqref{nonlexp} with \eqref{immderuslt}, we get \eqref{finderuslt} on an event $E$ of  the probability at least $1-2e^{-r}-p^{-A/2}$, 
\begin{align}\label{finderuslt}
\mathbf{Z}\leq L_0 \Big[ \Delta_{\beta}+\Delta_{-}+\Delta_{+}+\sqrt{\frac{r}{\log p}}\big(\Delta_{\beta}+\Delta_{-}\big)+\frac{ r}{\log p}\Delta_{+}\Big],\,\, \forall\,r\in(0, 2n),
\end{align}
where $L_0$ is some constant depending on $N_0$, $C_0$ and $\kappa$.

We will now choose $r=A\log p$ so as to obtain a weighted empirical process that holds uniformly over the constrains \eqref{splitnum},
\begin{align}\label{splitnum}
e^{-p}\leq \Delta_{\beta}\leq e^p,\,e^{-p}\leq \Delta_{-}\leq e^p,\,e^{-p}\leq \Delta_{+}\leq e^p.
\end{align}
To achieve this purpose, if we choose
$$
\Delta_{\beta}^k=\Delta_{-}^k=\Delta_{+}^k:=2^{-k},\quad k=-p,-p+1...,p-1,p,
$$
and
\begin{align}\label{splisec}
\Gamma_\Delta^{k,l,h}:=\Big\{\gamma,\,\,&\frac{1}{2}\Delta_\beta^k\leq \sqrt{\frac{\log p}{n}}\|\bbeta-\bbeta^*\|_1\leq \Delta_\beta^k, \frac{1}{2}\Delta_{-}^l\leq \gamma_n\| g-g^*\|_{2}\leq \Delta_{-}^l,\nonumber\\
& \frac{1}{2}\Delta_{+}^h\leq \gamma_n^2\| g-g^*\|_{K}\leq \Delta_{+}^h\Big\},
\end{align}
based on \eqref{finderuslt}, \eqref{weigheruslt} holds over an event $F(\Delta_{\beta}^k,\Delta^l_{-},\Delta^h_{+})$ of the probability $\mathbb{P}\big(F(\Delta_{\beta}^k,\Delta^l_{-},\Delta^h_{+})\big)\geq 1-p^{A/2}$ for any triplet $(\Delta_{\beta}^{k},\Delta^l_{-},\Delta^h_{+})$ satisfying \eqref{splitnum} and \eqref{splisec}, 
\begin{align}\label{weigheruslt}
\mathbf{Z}\leq L_0 \Big[ \Delta_{\beta}^k+\Delta_{-}^l+\Delta_{+}^h+\sqrt{A}\big(\Delta^k_{\beta}+\Delta^l_{-}\big)+A\Delta^h_{+}\Big],\,\, \forall\,r\in(0, 2n).
\end{align}
 On the event 
$E':=E\cap \big(\bigcap_{k,l,h}F(\Delta_{\beta}^k,\Delta^l_{-},\Delta^h_{+})\big)$, the intersection $\bigcap_{k,l,h}F(\Delta_{\beta}^k,\Delta^l_{-},\Delta^h_{+})$ is bounded by $(2p+1)^3$. Therefore, the lower bound of the probability $E'$ can be formulated as \eqref{P17},
\begin{align}
\label{P17}
\mathbb{P}(E')\geq 1-\mathbb{P}(E^c)- (2p+1)^3\exp(-A\log p)\geq 1-2p^{-A/2}-16p^{3-A}.
\end{align}
Thus, for any $k,l,h$, \eqref{P18} holds on the event  $E'$ with the construction of the function sets  $\Gamma_\Delta^{k,l,h}$ and \eqref{weigheruslt}, 
\begin{align}
\label{P18}
\mathbf{Z}\leq2L_0(1+\sqrt{A})\Big(\sqrt{\frac{\log p}{n}}\|\bbeta-\bbeta^*\|_1+\gamma_n\| g-g^*\|_{2}\Big)+
2L_0(1+A)\gamma_n^2\| g-g^*\|_{K}.
\end{align}
If it is true for either of the conditions $\Delta_{\beta}\leq e^{-p}$, or $\Delta_{-}\leq e^{-p}$ or $\Delta_{+}\leq e^{-p}$, it follows that \eqref{P18} with almost the same probability by monotonnicity of the left-hand side,
$$
\mathbf{Z}\leq2L_0(1+\sqrt{A})\Big(\sqrt{\frac{\log p}{n}}\|\bbeta-\bbeta^*\|_1+\gamma_n\| g-g^*\|_{2}\Big)+
2L_0(1+A)\gamma_n^2\| g-g^*\|_{K}+3e^{-p}.
$$
This completes the proof of Lemma \ref{wieghtineq}. \hspace*{10cm}$\Box$

\section{Discussion and Future Work}
In this paper, we have studied the estimation in the partially linear sparse  models for support vector machine, where the covariates split the linear component and the nonlinear component   within a reproducing kernel Hilbert space. 
An important feature of our analysis is that we develop a new weighted empirical process in the high dimensional semi-parametric setting, so that our derived rates are sharp  in comparison with the existing related results, even are comparable to  the least square estimations in the high dimensional setting.

There are some further research topics of this work. It is known that the parametric estimation for the partially linear mean regression does not depend on the functional complexity under some additional conditions.
This paper has not achieved such a better rate for the partially linear SVM. Therefore, how to improve the parametric estimation error is an interesting research problem. Besides, the lower bound for the partially linear SVM in the high dimensional setting have not been established under the margin information, which is an important complementary to our upper bounds.

\bigskip

\noindent\textbf{\Large Acknowledgement}\\
Yifan Xia's research is partially supported by
the National Social Science Fund of China (NSSFC-16BTJ013)(NSSFC-16ZDA010) and the Sichuan Project of Science \& Technology (2016JY0273). Shaogao Lv's research is partially supported by NSFC-11871277 and NSFC-11829101.

\bigskip
\bigskip


\end{document}